\documentclass[12pt]{article}
\author{Edward Hanson}
\title{\bf A characterization of bipartite Leonard\\ pairs using the notion of a tail}
\date{}
\usepackage{amssymb}
\usepackage{amsmath}
\usepackage{cite}

\newtheorem{definition}{Definition}[section]
\newtheorem{theorem}[definition]{Theorem}
\newtheorem{proposition}[definition]{Proposition}
\newtheorem{lemma}[definition]{Lemma}
\newtheorem{corollary}[definition]{Corollary}

\newtheorem{notation}[definition]{Notation}

\newtheorem{note}[definition]{Note}
\newtheorem{assumption}[definition]{Assumption}

\usepackage[margin=1.0in]{geometry}

\def\fld{\mathbb K}

\begin{document}
\maketitle

\begin{abstract}
\noindent Let $V$ denote a vector space with finite positive dimension. We consider an ordered pair of linear transformations 
$A: V\rightarrow V$ and $A^{*}: V\rightarrow V$ that satisfy (i) and (ii) below.
\begin{enumerate}
\item There exists a basis for $V$ with respect to which the matrix representing $A$ is irreducible tridiagonal and the matrix representing $A^{*}$ is diagonal.
\item There exists a basis for $V$ with respect to which the matrix representing $A^{*}$ is irreducible tridiagonal and the matrix representing $A$ is diagonal.
\end{enumerate}
We call such a pair a {\it Leonard pair} on $V$. Very roughly speaking, a Leonard pair is a linear algebraic abstraction of a $Q$-polynomial distance-regular graph. There is a well-known class of distance-regular graphs said to be bipartite and there is a related notion of a bipartite Leonard pair. Recently, M. S. Lang introduced the notion of a tail for bipartite distance-regular graphs, and there is an abstract version of this tail notion. Lang characterized the bipartite $Q$-polynomial distance-regular graphs using tails. In this paper, we obtain a similar characterization of the bipartite Leonard pairs using tails. Whereas Lang's arguments relied on the combinatorics of a distance-regular graph, our results are purely algebraic in nature.

\bigskip

\end{abstract}

\section{Introduction} \label{sec:intro}

We begin by recalling the notion of a Leonard pair \cite{T:subconst1, T:Leonard}. We will use the following terms. Let $X$ denote a square matrix. Then $X$ is called {\it tridiagonal} whenever each nonzero entry lies on either the diagonal, the subdiagonal, or the superdiagonal. Assume $X$ is tridiagonal. Then $X$ is called {\it irreducible} whenever each entry on
the subdiagonal is nonzero and each entry on the superdiagonal is nonzero.

\medskip

\noindent We now define a Leonard pair. For the rest of this paper, $\fld$ will denote a field.

\begin{definition} \label{def:lp} \rm \cite[Definition 1.1]{T:Leonard}
Let $V$ denote a vector space over $\fld$ with finite positive dimension. By a {\it Leonard pair} on $V$, we mean an ordered pair of linear transformations $A: V\rightarrow V$ and $A^{*}: V\rightarrow V$ that satisfy (i) and (ii) below.
\begin{enumerate}
\item There exists a basis for $V$ with respect to which the matrix representing $A$ is irreducible tridiagonal and the matrix representing $A^{*}$ is diagonal.
\item There exists a basis for $V$ with respect to which the matrix representing $A^{*}$ is irreducible tridiagonal and the matrix representing $A$ is diagonal.
\end{enumerate}
\end{definition}

\begin{note}
\rm
It is a common notational convention to use $A^{*}$ to represent the conjugate-transpose of $A$. We are not using this convention. In a Leonard pair $A,A^{*}$, the linear transformations $A$ and $A^{*}$ are arbitrary subject to (i), (ii) above.
\end{note}

\noindent Very roughly speaking, a Leonard pair is a linear algebraic abstraction of a $Q$-polynomial distance-regular graph \cite[p.~260]{BIbook} \cite[Definition 2.3]{T:subconst1}. In the theory of distance-regular graphs, there is a well-known set of parameters $a_{i}$ called intersection numbers. There also exists an abstract version of the $a_{i}$. A distance-regular graph is said to be bipartite whenever each $a_{i}$ equals zero. A bipartite Leonard pair is similarly defined.

\medskip

\noindent In \cite{Hanson, Hanson2}, we extended existing characterizations of $Q$-polynomial distance-regular graphs to obtain characterizations of Leonard pairs. We mention some details. In \cite{Lang}, M. S. Lang introduced the notion of a tail for bipartite distance-regular graphs. In \cite{JTZ}, the tail notion was applied to general distance-regular graphs. In \cite[Theorem 1.1]{JTZ}, these tails were used to characterize $Q$-polynomial distance-regular graphs. In \cite[Definition 4.5]{Hanson}, we introduced an abstract version of the tail notion and in \cite[Theorem 5.1]{Hanson}, we used it to characterize Leonard pairs. In \cite[Theorem 1.2]{Pascasio}, the $a_{i}$ were used to characterize $Q$-polynomial distance-regular graphs. In \cite{Hanson2}, we used the $a_{i}$ to characterize Leonard pairs; our main result was \cite[Theorem~5.1]{Hanson2}.

\medskip

\noindent We now summarize the present paper. Our point of departure is Lang's work concerning tails and bipartite distance-regular graphs \cite{Lang}. In \cite[Theorem 1.1]{Lang}, Lang used tails to characterize bipartite $Q$-polynomial distance-regular graphs. We use our abstract version of the tail notion \cite[Definition 4.5]{Hanson} to characterize bipartite Leonard pairs. Whereas Lang's arguments relied on the combinatorics of a distance-regular graph, our results are purely algebraic in nature.

\medskip

\noindent Our basic approach is as follows. We consider two linear transformations $A: V\rightarrow V$ and $A^*: V\rightarrow V$ that satisfy Definition \ref{def:lp}(i). We associate with $A,A^{*}$ a diagram $\Delta$. The vertices of $\Delta$ represent the eigenspaces of $A$ and the edges of $\Delta$ describe the action of $A^{*}$ on those eigenspaces. If $A,A^{*}$ is a Leonard pair, then $\Delta$ is a path. A tail in $\Delta$ is an ordered pair of distinct vertices $(i,j)$ such that $i$ is adjacent to no vertex in $\Delta$ besides $j$ and $j$ is adjacent to at most one vertex in $\Delta$ besides $i$. In our main result, we show that under mild assumptions, if $A$ is bipartite and $\Delta$ has a tail, then $A,A^{*}$ is a Leonard pair. In the proof, our primary step is to show that the eigenvalues of $A^{*}$ satisfy a three-term recurrence. This fact, together with \cite[Theorem 5.1]{Hanson}, implies that $\Delta$ is a path and this quickly yields that $A,A^{*}$ is a Leonard pair. Our main result is Theorem~\ref{thm:main}. 

\section{Leonard systems} \label{sec:LS}

When working with a Leonard pair, it is often convenient to consider a closely related object called a Leonard system. To prepare for our definition of a Leonard system, we recall a few concepts from linear algebra. From now on, we fix a nonnegative integer $d$. Let $\mbox{Mat}_{d+1}(\fld)$ denote the $\fld$-algebra consisting of all $d+1$ by $d+1$ matrices with entries in $\fld$. We index the rows and columns by $0,1,\ldots ,d$. Let $\fld^{d+1}$ denote the $\fld$-vector space consisting of all $d+1$ by $1$ matrices with entries in $\fld$. We index the rows by $0,1,\ldots ,d$. Recall that $\mbox{Mat}_{d+1}(\fld)$ acts on $\fld^{d+1}$ by left multiplication. Let $V$ denote a vector space over $\fld$ with dimension $d+1$. Let $\mbox{End}(V)$ denote the $\fld$-algebra consisting of all linear transformations from $V$ to $V$. For convenience, we abbreviate $\mathcal{A}=\mbox{End}(V)$. Observe that $\mathcal{A}$ is $\fld$-algebra isomorphic to $\mbox{Mat}_{d+1}(\fld)$ and that $V$ is irreducible as an $\mathcal{A}$-module. The identity of $\mathcal{A}$ will be denoted by $I$. Let $\{v_{i}\}_{i=0}^{d}$ denote a basis for $V$. For $X\in \mathcal{A}$ and $Y\in \mbox{Mat}_{d+1}(\fld)$, we say that $Y$ {\it represents} $X$ {\it with respect to} $\{v_{i}\}_{i=0}^{d}$ whenever $Xv_{j}=\sum_{i=0}^{d}Y_{ij}v_{i}$ for $0\leq j\leq d$. Let $A$ denote an element of $\mathcal{A}$. A subspace $W\subseteq V$ will be called an {\it eigenspace} of $A$ whenever $W\neq 0$ and there exists $\theta \in \fld$ such that $W=\{v\in V|Av=\theta v\}$; in this case, $\theta$ is the {\it eigenvalue} of $A$ associated with $W$. We say that $A$ is {\it diagonalizable} whenever $V$ is spanned by the eigenspaces of $A$. We say that $A$ is {\it multiplicity-free} whenever it has $d+1$ mutually distinct eigenvalues in $\fld$. Note that if $A$ is multiplicity-free, then $A$ is diagonalizable.

\begin{definition} \label{def:MOidem}
\rm
By a {\it system of mutually orthogonal idempotents} in $\mathcal{A}$, we mean a sequence $\{E_{i}\}_{i=0}^{d}$ of elements in $\mathcal{A}$ such that
\begin{equation}
E_{i}E_{j}=\delta_{i,j}E_{i}\qquad \qquad (0 \leq i,j \leq d), \notag
\end{equation}
\begin{equation}
{\rm rank}(E_{i})=1\qquad \qquad (0 \leq i \leq d). \notag
\end{equation}
\end{definition}

\begin{definition} \label{def:decomp}
\rm
By a {\it decomposition of $V$}, we mean a sequence $\{U_{i}\}_{i=0}^{d}$ consisting of one-dimensional subspaces of $V$ such that
\begin{equation}
V=\sum_{i=0}^{d}U_{i}\qquad \qquad \text{(direct sum)}. \notag
\end{equation}
\end{definition}

\noindent The following lemmas are routinely verified.

\begin{lemma} \label{lem:MOidem2}
Let $\{U_{i}\}_{i=0}^{d}$ denote a decomposition of $V$. For $0\leq i\leq d$, define $E_{i}\in \mathcal{A}$ such that $(E_{i}-I)U_{i}=0$ and $E_{i}U_{j}=0$ if $j\ne i$ $(0\leq j\leq d)$. Then $\{E_{i}\}_{i=0}^{d}$ is a system of mutually orthogonal idempotents. Conversely, given a system of mutually orthogonal idempotents $\{E_{i}\}_{i=0}^{d}$ in $\mathcal{A}$, define $U_{i}=E_{i}V$ for $0\leq i\leq d$. Then $\{U_{i}\}_{i=0}^{d}$ is a decomposition of $V$.
\end{lemma}

\begin{lemma} \label{lem:EsumI}
Let $\{E_{i}\}_{i=0}^{d}$ denote a system of mutually orthogonal idempotents in $\mathcal{A}$. Then $I=\sum_{i=0}^{d}E_{i}$.
\end{lemma}

\noindent Let $A$ denote a multiplicity-free element of $\mathcal{A}$ and let $\{\theta_{i}\}^{d}_{i=0}$ denote an ordering of the eigenvalues of $A$. For $0\leq i\leq d$, let $U_{i}$ denote the eigenspace of $A$ for $\theta_{i}$. Then $\{U_{i}\}_{i=0}^{d}$ is a decomposition of $V$; let $\{E_{i}\}_{i=0}^{d}$ denote the corresponding system of idempotents from Lemma \ref{lem:MOidem2}. One checks that $A=\sum_{i=0}^{d}\theta_{i}E_{i}$ and $AE_{i}=E_{i}A=\theta_{i}E_{i}$ for $0\leq i\leq d$. Moreover,
\begin{equation}
E_{i}=\prod_{\genfrac{}{}{0pt}{}{0 \leq  j \leq d}{j\not=i}}\frac{A-\theta_{j}I}{\theta_{i}-\theta_{j}}\qquad \qquad (0\leq i\leq d). \notag
\end{equation}
We refer to $E_{i}$ as the {\it primitive idempotent} of $A$ corresponding to $U_{i}$ (or $\theta_{i}$).

\medskip

\noindent We now define a Leonard system.

\begin{definition} \label{def:ls} \rm \cite[Definition 1.4]{T:Leonard}
By a {\it Leonard system} on $V$, we mean a sequence 
\begin{equation}
(A; \{E_{i}\}_{i=0}^{d}; A^{*}; \{E^{*}_{i}\}_{i=0}^{d}) \notag
\end{equation}
which satisfies (i)--(v) below.
\begin{enumerate}
\item Each of $A,A^{*}$ is a multiplicity-free element of $\mathcal{A}$.
\item $\{E_{i}\}_{i=0}^{d}$ is an ordering of the primitive idempotents of $A$.
\item $\{E^{*}_{i}\}_{i=0}^{d}$ is an ordering of the primitive idempotents of $A^{*}$.
\item ${\displaystyle{E^{*}_{i}AE^{*}_{j} =
\begin{cases}
0, & \text{if $\;|i-j|>1$;} \\
\neq 0, & \text{if $\;|i-j|=1$}
\end{cases}
}}
\qquad \qquad (0 \leq i,j\leq d).$
\item ${\displaystyle{E_{i}A^{*}E_{j} =
\begin{cases}
0, & \text{if $\;|i-j|>1$;} \\
\neq 0, & \text{if $\;|i-j|=1$}
\end{cases}
}}
\qquad \qquad (0 \leq i,j\leq d).$
\end{enumerate}
\end{definition}

\noindent Leonard systems and Leonard pairs are related as follows. Let $(A; \{E_{i}\}_{i=0}^{d}; A^{*}; \{E^{*}_{i}\}_{i=0}^{d})$ denote a Leonard system on $V$. For $0\leq i\leq d$, let $v_{i}$ denote a nonzero vector in $E_{i}V$. Then the sequence $\{v_{i}\}_{i=0}^{d}$ is a basis for $V$ which satisfies Definition \ref{def:lp}(ii). For $0\leq i\leq d$, let $v^{*}_{i}$ denote a nonzero vector in $E^{*}_{i}V$. Then the sequence $\{v^{*}_{i}\}_{i=0}^{d}$ is a basis for $V$ which satisfies Definition \ref{def:lp}(i). By these comments, the pair $A,A^{*}$ is a Leonard pair on $V$. Conversely, let $A,A^{*}$ denote a Leonard pair on $V$. By \cite[Lemma 1.3]{T:Leonard}, each of $A,A^{*}$ is multiplicity-free. Let $\{v_{i}\}_{i=0}^{d}$ denote a basis for $V$ which satisfies Definition \ref{def:lp}(ii). For $0\leq i\leq d$, the vector $v_{i}$ is an eigenvector for $A$; let $E_{i}$ denote the corresponding primitive idempotent. Let $\{v^{*}_{i}\}_{i=0}^{d}$ denote a basis for $V$ which satisfies Definition \ref{def:lp}(i). For $0\leq i\leq d$, the vector $v^{*}_{i}$ is an eigenvector for $A^{*}$; let $E^{*}_{i}$ denote the corresponding primitive idempotent. Then $(A; \{E_{i}\}_{i=0}^{d}; A^{*}; \{E^{*}_{i}\}_{i=0}^{d})$ is a Leonard system on $V$.

\medskip

\noindent We make some observations. Let $(A; \{E_{i}\}_{i=0}^{d}; A^{*}; \{E^{*}_{i}\}_{i=0}^{d})$ denote a Leonard system on $V$. For $0\leq i\leq d$, let $\theta_{i}$ (resp. $\theta^{*}_{i}$) denote the eigenvalue of $A$ (resp. $A^{*}$) associated with $E_{i}V$ (resp. $E^{*}_{i}V$). By construction, $\{\theta_{i}\}_{i=0}^{d}$ (resp. $\{\theta^{*}_{i}\}_{i=0}^{d}$) are mutually distinct and contained in $\fld$. By \cite[Theorem 12.7]{T:Leonard}, the expressions
\begin{equation} \label{eq:thetarecur}
\frac{\theta_{i-2}-\theta_{i+1}}{\theta_{i-1}-\theta_{i}}, \qquad \frac{\theta^{*}_{i-2}-\theta^{*}_{i+1}}{\theta^{*}_{i-1}-\theta^{*}_{i}}
\end{equation}
are equal and independent of $i$ for $2\leq i\leq d-1$.

\section{Basic assumptions}

In this section, we establish the basic setting in which we work.

\begin{assumption} \label{assum:E*AE*}
\rm
Let $\{E^{*}_{i}\}^{d}_{i=0}$ denote a system of mutually orthogonal idempotents in $\mathcal{A}$. Let $A$ denote an element of $\mathcal{A}$ such that
\begin{equation} \label{eq:E*AE*}
{\displaystyle{E^{*}_{i}AE^{*}_{j} =
\begin{cases}
0, & \text{if $\;|i-j|>1$;} \\
\neq 0, & \text{if $\;|i-j|=1$}
\end{cases}
}}
\qquad \qquad (0 \leq i,j\leq d).
\end{equation}
\end{assumption}

\begin{definition} \label{def:a}
\rm
With reference to Assumption \ref{assum:E*AE*}, define
\begin{equation}
a_{i}=\mbox{tr}(E^{*}_{i}A) \qquad \qquad (0 \leq i \leq d), \notag
\end{equation}
where $\mbox{tr}$ denotes trace.
\end{definition}

\begin{proposition} \label{prop:E*AE*} {\rm \cite[Proposition 3.6]{Hanson2}}
With reference to Assumption \ref{assum:E*AE*}, $E^{*}_{i}AE^{*}_{i}=a_{i}E^{*}_{i}$ for $0\leq i\leq d$.
\end{proposition}

\noindent We have been discussing the situation of Assumption \ref{assum:E*AE*}. We now modify this situation as follows.

\begin{assumption} \label{assum:A*}
\rm
Let $A$ and $\{E^{*}_{i}\}_{i=0}^{d}$ be as in Assumption \ref{assum:E*AE*}. Furthermore, assume that $A$ is multiplicity-free, with primitive idempotents $\{E_{i}\}_{i=0}^{d}$ and eigenvalues $\{\theta_{i}\}_{i=0}^{d}$. Additionally, let $\{\theta^{*}_{i}\}^{d}_{i=0}$ denote scalars in $\fld$ and let $A^{*}=\sum_{i=0}^{d}\theta^{*}_{i}E^{*}_{i}$. To avoid trivialities, assume that $d\geq 1$.
\end{assumption}

\begin{lemma} \label{lem:undirected} {\rm \cite[Lemma 3.8]{Hanson}}
With reference to Assumption \ref{assum:A*} and for $0\leq i,j\leq d$, $E_{i}A^{*}E_{j}=0$ if and only if $E_{j}A^{*}E_{i}=0$.
\end{lemma}

\section{The graph $\Delta$} \label{sec:delta}

In the following discussion, a graph is understood to be finite and undirected, without loops or multiple edges.

\begin{definition} \label{def:delta}
\rm
With reference to Assumption \ref{assum:A*}, let $\Delta$ denote the graph with vertex set $\{0, 1,\ldots ,d\}$ such that two vertices $i, j$ are adjacent if and only if $i\neq j$ and $E_{i}A^{*}E_{j}\neq 0$. The graph $\Delta$ is well-defined in view of Lemma \ref{lem:undirected}.
\end{definition}

\begin{lemma} \label{lem:lspath} {\rm \cite[Lemma 4.2]{Hanson}}
With reference to Assumption \ref{assum:A*}, the following {\rm (i)}, {\rm (ii)} are equivalent.
\begin{enumerate}
\item[\rm (i)] The sequence $(A; \{E_{i}\}^{d}_{i=0}; A^{*}; \{E^{*}_{i}\}^{d}_{i=0})$ is a Leonard system.
\item[\rm (ii)] The diagram $\Delta$ is a path such that vertices $i-1, i$ are adjacent for $1\leq i\leq d$.
\end{enumerate}
\end{lemma}

\begin{definition} \label{def:Qpoly1}
\rm
With reference to Assumption \ref{assum:A*}, the given ordering $\{E_{i}\}_{i=0}^{d}$ of the primitive idempotents of $A$ is said to be {\it $Q$-polynomial} whenever the equivalent conditions (i), (ii) hold in Lemma \ref{lem:lspath}.
\end{definition}

\begin{definition} \label{def:Qpoly2}
\rm
With reference to Assumption \ref{assum:A*}, let $(E,F)$ denote an ordered pair of distinct primitive idempotents for $A$. This pair will be called {\it $Q$-polynomial} whenever there exists a $Q$-polynomial ordering $\{E_{i}\}_{i=0}^{d}$ of the primitive idempotents of $A$ such that $E=E_{0}$ and $F=E_{1}$.
\end{definition}

\noindent The following is motivated by \cite[Definition 5.1]{Lang}.

\begin{definition} \label{def:tail}
\rm
With reference to Assumption \ref{assum:A*}, let $(E,F)=(E_{i},E_{j})$ denote an ordered pair of distinct primitive idempotents for $A$. This pair will be called a {\it tail} whenever the following (i), (ii) occurs in $\Delta$.
\begin{enumerate}
\item[\rm (i)] $i$ is adjacent to no vertex in $\Delta$ besides $j$.
\item[\rm (ii)] $j$ is adjacent to at most one vertex in $\Delta$ besides $i$.
\end{enumerate}
\end{definition}

\begin{lemma} \label{lem:Qpolytail}
With reference to Assumption \ref{assum:A*}, let $(E,F)$ denote an ordered pair of distinct primitive idempotents for $A$. If $(E,F)$ is $Q$-polynomial, then $(E,F)$ is a tail.
\end{lemma}
\noindent {\it Proof:} Compare Definitions \ref{def:Qpoly1} and \ref{def:tail}. \hfill $\Box$ \\

\noindent When working with a tail, pertinent information can be obtained by considering the following related notion.

\begin{definition} \label{def:stail}
\rm
With reference to Assumption \ref{assum:A*}, let $E=E_{i}$ denote a primitive idempotent for $A$. This idempotent will be called a {\it leaf} whenever $i$ is adjacent to at most one vertex in $\Delta$.
\end{definition}

\noindent With reference to Assumption \ref{assum:A*}, let $(E,F)$ denote an ordered pair of distinct primitive idempotents for $A$. Note that by Definitions \ref{def:tail} and \ref{def:stail}, if $(E,F)$ is a tail then $E$ is a leaf.

\medskip

\noindent We now discuss the connectivity of $\Delta$.

\begin{proposition} \label{prop:deltaconnected} {\rm \cite[Proposition 4.9]{Hanson2}}
With reference to Assumption \ref{assum:A*}, assume further that $\theta^{*}_{i}\neq \theta^{*}_{0}$ for $1\leq i \leq d$. Then $\Delta$ is connected.
\end{proposition}

\noindent Recall the following theorem that characterizes the $Q$-polynomial property using tails.
 
\begin{theorem} {\rm \cite[Theorem 5.1]{Hanson}} \label{thm:tail}
With reference to Assumption \ref{assum:A*}, let $(E,F)$ denote an ordered pair of distinct primitive idempotents for A. Then this pair is $Q$-polynomial if and only if the following {\rm (i)--(iii)} hold.
\begin{enumerate}
\item[\rm (i)] $(E,F)$ is a tail.
\item[\rm (ii)] There exists $\beta \in \fld$ such that $\theta^{*}_{i-1}-\beta \theta^{*}_{i}+\theta^{*}_{i+1}$ is independent of $i$ for $1\leq i \leq d-1$.
\item[\rm (iii)] $\theta^{*}_{i}\neq \theta^{*}_{0}$ for $1\leq i \leq d$.
\end{enumerate}
\end{theorem}

\noindent For later use, we now review some results from \cite{Hanson2}. We first make a necessary definition.

\begin{definition} \label{def:a*}
\rm
With reference to Assumption \ref{assum:A*}, define
\begin{equation}
a^{*}_{i}=\mbox{tr}(E_{i}A^{*}) \qquad \qquad (0 \leq i \leq d). \notag
\end{equation}
\end{definition}

\begin{proposition} \label{prop:EA*E} {\rm \cite[Proposition 6.2]{Hanson2}}
With reference to Assumption \ref{assum:A*}, $E_{i}A^{*}E_{i}=a^{*}_{i}E_{i}$ for $0\leq i\leq d$.
\end{proposition}

\begin{lemma} \label{lem:A*-a*I} {\rm \cite[Lemma 6.3]{Hanson2}}
With reference to Assumption \ref{assum:A*}, the following {\rm (i)}, {\rm (ii)} are equivalent.
\begin{enumerate}
\item[\rm (i)] In the diagram $\Delta$, vertex $0$ is adjacent to vertex $1$ and no other vertices.
\item[\rm (ii)] There exists $\kappa\in\fld$ such that $(A^{*}-\kappa I)E_{0}V=E_{1}V$.
\end{enumerate}
Suppose conditions {\rm (i)} and {\rm (ii)} hold. Then $\kappa=a^{*}_{0}$.
\end{lemma}

\section{The bipartite case}

\begin{definition} \label{def:bipartite}
\rm
With reference to Assumption \ref{assum:E*AE*}, we say that $A$ is {\it bipartite} whenever $a_{i}=0$ for $0\leq i\leq d$.
\end{definition}

\begin{assumption} \label{assum:A_bip}
\rm
Let $A$, $A^{*}$, $\{E_{i}\}_{i=0}^{d}$, and $\{E^{*}_{i}\}_{i=0}^{d}$ be as in Assumption \ref{assum:A*}. Furthermore, assume that $A$ is bipartite.
\end{assumption}

\noindent We will need material from \cite{Hanson2} in a form appropriate to the case when $A$ is bipartite. Therefore, the relevant results will be presented with reference to Assumption \ref{assum:A_bip}.

\begin{definition} \label{def:feasible}
\rm
With reference to Assumption \ref{assum:A_bip}, let $\{v_{i}\}_{i=0}^{d}$ denote a basis of $V$. We say that this basis is {\it feasible} whenever $v_{i}\in E^{*}_{i}V$ for $0\leq i\leq d$.
\end{definition}

\noindent With reference to Assumption \ref{assum:A_bip}, let $\{v_{i}\}_{i=0}^{d}$ denote a feasible basis of $V$. By Definition~\ref{def:bipartite}, Proposition \ref{prop:E*AE*}, and Definition \ref{def:feasible}, the matrices representing $A$ and $A^{*}$ with respect to $\{v_{i}\}_{i=0}^{d}$ are
\begin{equation} \label{eq:A_A*_mat_gen}
A:\left(
\begin{array}
{ c c c c c c}
  0     & b_{0} &       &       &       & {\bf 0} \\
  c_{1} & 0     & b_{1} &       &       & \\
        & c_{2} & \cdot & \cdot &       & \\
        &       & \cdot & \cdot & \cdot & \\
        &       &       & \cdot & \cdot & b_{d-1} \\
{\bf 0} &       &       &       & c_{d} & 0
\end{array}
\right)
\qquad
A^{*}:\left(
\begin{array}
{ c c c c c c}
\theta^{*}_{0} &                &       &       &       & {\bf 0} \\
               & \theta^{*}_{1} &       &       &       & \\
               &                & \cdot &       &       & \\
               &                &       & \cdot &       & \\
               &                &       &       & \cdot & \\
{\bf 0}        &                &       &       &       & \theta^{*}_{d}
\end{array}
\right),
\end{equation}
where each of the scalars $\{b_{i}\}_{i=0}^{d-1}$ and $\{c_{i}\}_{i=1}^{d}$ is nonzero. For notational convenience, let $b_{d}=0$ and $c_{0}=0$.

\medskip

\noindent Observe that for $1\leq i\leq d$, the product $b_{i-1}c_{i}$ is independent of our choice of feasible basis. However, the scalars $\{b_{i}\}_{i=0}^{d-1}$ depend on the choice of feasible basis, as shown in the following lemma.

\begin{lemma} \label{lem:A_A*_mat2} {\rm \cite[Lemma 7.2]{Hanson2}}
With reference to Assumption \ref{assum:A_bip}, let $\{\beta_{i}\}_{i=0}^{d-1}$ denote a sequence of nonzero scalars taken from $\fld$. Then there exists a feasible basis of $V$ with respect to which the matrix representing $A$ has $(i,i+1)$-entry $\beta_{i}$ for $0\leq i\leq d-1$.
\end{lemma}

\begin{definition}
\rm
With reference to Assumption \ref{assum:A_bip}, it follows by Lemma \ref{lem:A_A*_mat2} that there exists a feasible basis of $V$ such that $b_{i}=1$ for $0\leq i\leq d-1$. We call this basis the \emph{normalized feasible basis of $V$}.
\end{definition}

\section{Normalizing idempotents}

Let $\lambda$ denote an indeterminate. Let $\fld[\lambda]$ denote the $\fld$-algebra consisting of the polynomials in $\lambda$ that have all coefficients in $\fld$.

\begin{definition} \label{def:u_i}
\rm
With reference to Assumption \ref{assum:A_bip}, let $\{v_{i}\}_{i=0}^{d}$ denote a feasible basis for $V$. Define a sequence of polynomials $\{u_{i}\}_{i=0}^{d+1}$ in $\fld[\lambda]$ by
\begin{align}
u_{0} &= 1, \label{eq:u_0} \\
\lambda u_{i} &= c_{i}u_{i-1}+b_{i}u_{i+1} \qquad \qquad (0 \leq i \leq d-1), \label{eq:u_i} \\
\lambda u_{d} &= c_{d}u_{d-1}+\frac{u_{d+1}}{b_{0}b_{1}\cdots b_{d-1}}, \label{eq:u_d+1}
\end{align}
where $u_{-1}=0$. Observe that for $0\leq i\leq d+1$, the polynomial $u_{i}$ has degree $i$. Moreover, the coefficient of $\lambda^{i}$ in $u_{i}$ equals $(b_{0}b_{1}\cdots b_{i-1})^{-1}$ if $0\leq i\leq d$ and $1$ if $i=d+1$. We say that the sequence $\{u_{i}\}_{i=0}^{d+1}$ \emph{corresponds} to the feasible basis $\{v_{i}\}_{i=0}^{d}$. 
\end{definition}

\begin{definition} \label{def:p_i}
\rm
With reference to Assumption \ref{assum:A_bip}, let $\{p_{i}\}_{i=0}^{d+1}$ denote the polynomial sequence that corresponds to the normalized feasible basis of $V$. Observe that $p_{i}$ is monic of degree $i$ for $0\leq i\leq d+1$.
\end{definition}

\noindent We adopt the following assumption.

\begin{assumption} \label{assum:fix_feas}
\rm
With reference to Assumption \ref{assum:A_bip}, fix a feasible basis $\{v_{i}\}_{i=0}^{d}$ of $V$. Let $\{b_{i}\}_{i=0}^{d-1}$ and $\{c_{i}\}_{i=1}^{d}$ denote the corresponding scalars from (\ref{eq:A_A*_mat_gen}). Let $\{u_{i}\}_{i=0}^{d+1}$ denote the corresponding sequence of polynomials from Definition \ref{def:u_i}.
\end{assumption}

\noindent Recall the polynomials $\{p_{i}\}_{i=0}^{d+1}$ from Definition \ref{def:p_i}. From the perspective of Assumption \ref{assum:fix_feas}, these polynomials appear as follows.

\begin{lemma} {\rm \cite[Lemma 7.7]{Hanson2}}
With reference to Assumptions \ref{assum:A_bip} and \ref{assum:fix_feas},
\begin{align}
p_{0} &= 1, \label{eq:p_0} \\
\lambda p_{i} &= b_{i-1}c_{i}p_{i-1}+p_{i+1} \qquad \qquad (0 \leq i \leq d), \label{eq:p_i}
\end{align}
where $p_{-1}=0$.
\end{lemma}

\begin{lemma} \label{lem:p_u} {\rm \cite[Lemma 7.8]{Hanson2}}
With reference to Assumptions \ref{assum:A_bip} and \ref{assum:fix_feas},
\begin{align}
  u_{i} &= \frac{p_{i}}{b_{0}b_{1}\cdots b_{i-1}} \qquad \qquad (0\leq i\leq d), \notag \\
u_{d+1} &= p_{d+1}. \notag
\end{align}
\end{lemma}
\noindent {\it Proof:} Compare (\ref{eq:u_0})--(\ref{eq:u_d+1}) with (\ref{eq:p_0}) and (\ref{eq:p_i}). \hfill $\Box$ \\

\begin{lemma} \label{lem:recur_u} {\rm \cite[Lemma 7.9]{Hanson2}}
With reference to Assumptions \ref{assum:A_bip} and \ref{assum:fix_feas}, let $v$ denote a nonzero vector in $V$ and write $v=\sum_{i=0}^{d}\alpha_{i}v_{i}$. Let $\theta \in \fld$. Then the following {\rm (i)--(iii)} are equivalent.
\begin{enumerate}
\item[\rm (i)] The vector $v$ is an eigenvector for $A$ with eigenvalue $\theta$.
\item[\rm (ii)] For $0\leq i\leq d$,
\begin{equation}
c_{i}\alpha_{i-1}+b_{i}\alpha_{i+1}=\theta\alpha_{i}, \notag
\end{equation}
where $\alpha_{-1}$ and $\alpha_{d+1}$ are indeterminates.
\item[\rm (iii)] $\alpha_{i}=u_{i}(\theta)\alpha_{0}$ for $0\leq i\leq d$ and $u_{d+1}(\theta)=0$.
\end{enumerate}
Suppose conditions {\rm (i)--(iii)} hold. Then $\alpha_{0}\neq 0$.
\end{lemma}

\noindent With reference to Assumptions \ref{assum:A_bip} and \ref{assum:fix_feas}, note that by \cite[Corollary 7.10]{Hanson2}, the polynomial $u_{d+1}$ is the characteristic polynomial for $A$.

\medskip

\noindent With reference to Assumptions \ref{assum:A_bip} and \ref{assum:fix_feas}, let $\theta$ denote an eigenvalue of $A$. In Lemma \ref{lem:recur_u}(iii), we encountered the sequence $\{u_{i}(\theta)\}_{i=0}^{d}$. In the theory of distance-regular graphs, this sequence is called the cosine sequence for $\theta$. Motivated by this, we call the sequence $\{u_{i}(\theta)\}_{i=0}^{d}$ the \emph{cosine sequence for $\theta$ with respect to $\{v_{i}\}_{i=0}^{d}$}. Sometimes it is clear from the context what the basis $\{v_{i}\}_{i=0}^{d}$ is. In this case, we will refer to $\{u_{i}(\theta)\}_{i=0}^{d}$ as the \emph{cosine sequence for $\theta$}.

\begin{lemma} \label{lem:cos_recur} {\rm \cite[Lemma 7.11]{Hanson2}}
With reference to Assumptions \ref{assum:A_bip} and \ref{assum:fix_feas}, let $\theta$ and $\{\alpha_{i}\}_{i=0}^{d}$ denote scalars in $\fld$. Then the following {\rm (i)}, {\rm (ii)} are equivalent.
\begin{enumerate}
\item[\rm (i)] The scalar $\theta$ is an eigenvalue for $A$ and $\{\alpha_{i}\}_{i=0}^{d}$ is the corresponding cosine sequence.
\item[\rm (ii)] $\alpha_{0}=1$ and
\begin{equation} \label{eq:cos_TTR}
c_{i}\alpha_{i-1}+b_{i}\alpha_{i+1}=\theta\alpha_{i} \qquad (0 \leq i \leq d),
\end{equation}
where $\alpha_{-1}$ and $\alpha_{d+1}$ are indeterminates.
\end{enumerate}
Suppose conditions {\rm (i)} and {\rm (ii)} hold. Then $\sum_{i=0}^{d}\alpha_{i}v_{i}$ is an eigenvector for $A$ with eigenvalue $\theta$.
\end{lemma}
\noindent {\it Proof:} This is an immediate consequence of Lemma \ref{lem:recur_u}. \hfill $\Box$ \\

\begin{proposition} \label{prop:genstailTTR} {\rm \cite[Proposition 7.13]{Hanson2}}
With reference to Assumptions \ref{assum:A_bip} and \ref{assum:fix_feas}, the following {\rm (i)}, {\rm (ii)} are equivalent.
\begin{enumerate}
\item[\rm (i)] In the diagram $\Delta$, vertex $0$ is adjacent to vertex $1$ and no other vertices.
\item[\rm (ii)] The cosine sequence $\{\alpha_{i}\}_{i=0}^{d}$ for $\theta_{0}$ satisfies
\begin{equation} \label{eq:genstailTTR}
c_{i}\theta^{*}_{i-1}\alpha_{i-1}+b_{i}\theta^{*}_{i+1}\alpha_{i+1}-\theta_{0}\theta^{*}_{i}\alpha_{i}=(\theta_{1}-\theta_{0})(\theta^{*}_{i}-a^{*}_{0})\alpha_{i} \quad (0\leq i\leq d),
\end{equation}
where each of $\alpha_{-1}$, $\alpha_{d+1}$, $\theta^{*}_{-1}$, and $\theta^{*}_{d+1}$ is indeterminate. Furthermore, there exists an integer $i$ ($0\leq i\leq d$) such that the right-hand side of {\rm (\ref{eq:genstailTTR})} is not equal to $0$.
\end{enumerate}
\end{proposition}

\begin{proposition} \label{prop:rowsum2}
With reference to Assumptions \ref{assum:A_bip} and \ref{assum:fix_feas}, the following {\rm (i)--(iii)} are equivalent for $\theta \in \fld$.
\begin{enumerate}
\item[\rm (i)] There exists a feasible basis for $V$ with respect to which the matrix on the left in {\rm (\ref{eq:A_A*_mat_gen})} has constant row sum $\theta$.
\item[\rm (ii)] The scalar $\theta$ is an eigenvalue of $A$ and $u_{i}(\theta)\neq 0$ for $0\leq i\leq d$, where the polynomials $\{u_{i}\}_{i=0}^{d}$ are from Assumption \ref{assum:fix_feas}.
\item[\rm (iii)] The scalar $\theta$ is an eigenvalue of $A$ and $p_{i}(\theta)\neq 0$ for $0\leq i\leq d$, where the polynomials $\{p_{i}\}_{i=0}^{d}$ are from Definition \ref{def:p_i}.
\end{enumerate}
\end{proposition}
\noindent {\it Proof:} (i) $\Leftrightarrow$ (ii). This is \cite[Proposition 7.16]{Hanson2}.

\medskip

\noindent (ii) $\Leftrightarrow$ (iii). This is an immediate consequence of Lemma \ref{lem:p_u}. \hfill $\Box$ \\

\noindent We make a comment on Proposition \ref{prop:rowsum2}. Although condition (ii) seemingly depends on a choice of feasible basis, it is clear that the equivalent conditions (i) and (iii) do not. Therefore, Proposition \ref{prop:rowsum2} can be referenced without referencing Assumption \ref{assum:fix_feas}, as in the following definitions.

\begin{definition} \label{def:normalizing}
\rm
With reference to Assumption \ref{assum:A_bip}, suppose $\theta$ is an eigenvalue of $A$. We say that $\theta$ is {\it normalizing} whenever the equivalent conditions (i)--(iii) hold in Proposition \ref{prop:rowsum2}.
\end{definition}

\begin{definition} \label{def:normalizing2}
\rm
With reference to Assumption \ref{assum:A_bip}, let $E$ denote a primitive idempotent for $A$. This idempotent will be called {\it normalizing} whenever $E$ corresponds to a normalizing eigenvalue $\theta$.
\end{definition}

\section{Algebraic consequences of the leaf condition}

For the remainder of this paper, we fix the following notation.

\begin{notation} \label{not:alpha}
\rm
With reference to Assumptions \ref{assum:A_bip} and \ref{assum:fix_feas}, let $\{\alpha_{i}\}_{i=0}^{d}$ denote the cosine sequence for $\theta_{0}$. For notational convenience, let $\alpha_{-1}$ and $\alpha_{d+1}$ denote indeterminates.
\end{notation}

\noindent With reference to Assumption \ref{assum:A_bip}, Assumption \ref{assum:fix_feas}, and Notation \ref{not:alpha}, it follows by (\ref{eq:cos_TTR}) that
\begin{equation} \label{eq:cos_TTR2}
c_{i}\alpha_{i-1}+b_{i}\alpha_{i+1} = \theta_{0}\alpha_{i}.
\end{equation}
Letting $i=0$ and $i=d$ in (\ref{eq:cos_TTR2}), we obtain
\begin{align}
b_{0}\alpha_{1} &= \theta_{0}\alpha_{0}, \label{eq:b0alpha1} \\
c_{d}\alpha_{d-1} &= \theta_{0}\alpha_{d}. \label{eq:cdalphad-1}
\end{align}

\medskip

\noindent Now suppose that in the diagram $\Delta$, vertex $0$ is adjacent to vertex $1$ and no other vertices. By (\ref{eq:genstailTTR}), the following holds for $0\leq i\leq d$,
\begin{equation} \label{eq:genstailTTR2}
c_{i}(\theta^{*}_{i-1}-\theta^{*}_{0})\alpha_{i-1}+b_{i}(\theta^{*}_{i+1}-\theta^{*}_{0})\alpha_{i+1} = \theta_{1}(\theta^{*}_{i}-\theta^{*}_{0})\alpha_{i}+(\theta_{1}-\theta_{0})(\theta^{*}_{0}-a^{*}_{0})\alpha_{i}.
\end{equation}

\medskip

\noindent In the sections that follow, we discuss a method for obtaining the three-term recurrence relationship on $\{\theta^{*}_{i}\}_{i=0}^{d}$ referenced in Theorem \ref{thm:tail} as condition (ii). Note that if $d<3$, then this relationship is automatically satisfied, so we focus our attention on the case when $d\geq 3$.

\begin{assumption} \label{assum:A_bip_d>2}
\rm
Let $A$, $A^{*}$, $\{E_{i}\}_{i=0}^{d}$, and $\{E^{*}_{i}\}_{i=0}^{d}$ be as in Assumption \ref{assum:A_bip}. Furthermore, assume that $d\geq 3$.
\end{assumption}

\begin{note} \label{note:assumptions}
\rm
With reference to Assumption \ref{assum:A_bip_d>2}, Theorem \ref{thm:main} involves the following conditions (i)--(iii).
\begin{enumerate}
\item[(i)] The primitive idempotent $E$ is normalizing.
\item[(ii)] $(E,F)$ is a tail.
\item[(iii)] $\{\theta^{*}_{i}\}_{i=0}^{d}$ are mutually distinct.
\end{enumerate}
Throughout the remainder of this paper, we present results that depend on various combinations of these three conditions. Therefore, we will reference which of these conditions are necessary assumptions for each result. When assuming (i), we take $E=E_{0}$ without loss of generality. When assuming (ii), we take $(E,F)=(E_{0},E_{1})$ without loss of generality. Furthermore, when we assume both (ii) and (iii), the primitive idempotent $E_{1}$ is adjacent to one primitive idempotent in $\Delta$ besides $E_{0}$, which we denote $E_{2}$.
\end{note}

\noindent We make some additional comments on Note \ref{note:assumptions}. By Notation \ref{not:alpha}, $E_{0}$ is normalizing if and only if $\alpha_{i}\neq 0$ for $0\leq i\leq d$. Also, note that under the assumption of Note \ref{note:assumptions}(iii), the diagram $\Delta$ is connected by Proposition \ref{prop:deltaconnected}.

\medskip

\noindent We now proceed to derive expressions for the scalars $a^{*}_{0}$.

\begin{lemma}
With reference to Assumption \ref{assum:A_bip_d>2}, further assume condition {\rm (i)} from Note \ref{note:assumptions}. Suppose that in the diagram $\Delta$, vertex $0$ is adjacent to vertex $1$ and no other vertices. Then
\begin{equation}
a^{*}_{0} = \frac{\theta_{1}\theta^{*}_{0}-\theta_{0}\theta^{*}_{1}}{\theta_{1}-\theta_{0}} = \frac{\theta_{1}\theta^{*}_{d}-\theta_{0}\theta^{*}_{d-1}}{\theta_{1}-\theta_{0}}. \label{eq:a*0}
\end{equation}
\end{lemma}
\noindent {\it Proof:} Without loss of generality, fix a feasible basis as in Assumption \ref{assum:fix_feas}. Substituting (\ref{eq:b0alpha1}) into (\ref{eq:genstailTTR2}) at $i=0$ and using the fact that $\theta_{0}$ is normalizing by Note \ref{note:assumptions}(i), we obtain the first part of (\ref{eq:a*0}). Similarly, the second part of (\ref{eq:a*0}) is a consequence of (\ref{eq:genstailTTR2}) at $i=d$ and (\ref{eq:cdalphad-1}). \hfill $\Box$ \\

\begin{corollary} \label{cor:theta01}
With reference to Assumption \ref{assum:A_bip_d>2}, further assume conditions {\rm (i)} and {\rm (iii)} from Note \ref{note:assumptions}. Suppose that in the diagram $\Delta$, vertex $0$ is adjacent to vertex $1$ and no other vertices. Then
\begin{equation}
\theta_{0}(\theta^{*}_{d-1}-\theta^{*}_{1}) = \theta_{1}(\theta^{*}_{d}-\theta^{*}_{0}). \label{eq:theta01}
\end{equation}
Moreover, $\theta_{0}$ and $\theta_{1}$ are both nonzero.
\end{corollary}
\noindent {\it Proof:} First note that (\ref{eq:theta01}) is a consequence of (\ref{eq:a*0}). By Note \ref{note:assumptions}(iii), both of the coefficients of $\theta_{0}$ and $\theta_{1}$ in (\ref{eq:theta01}) are nonzero. Therefore, both $\theta_{0}$ and $\theta_{1}$ are nonzero. \hfill $\Box$ \\

\begin{lemma}
With reference to Assumptions \ref{assum:A_bip_d>2} and \ref{assum:fix_feas}, further assume conditions {\rm (i)} and {\rm (iii)} from Note \ref{note:assumptions}. Suppose that in the diagram $\Delta$, vertex $0$ is adjacent to vertex $1$ and no other vertices. Then
\begin{align} 
b_{i} &= \frac{\theta_{1}(\theta^{*}_{i}-\theta^{*}_{0})-\theta_{0}(\theta^{*}_{i-1}-\theta^{*}_{1})}{\theta^{*}_{i+1}-\theta^{*}_{i-1}}\frac{\alpha_{i}}{\alpha_{i+1}} \qquad \qquad (1\leq i\leq d-1), \label{eq:bi_1} \\
c_{i} &= \frac{\theta_{0}(\theta^{*}_{i+1}-\theta^{*}_{1})-\theta_{1}(\theta^{*}_{i}-\theta^{*}_{0})}{\theta^{*}_{i+1}-\theta^{*}_{i-1}}\frac{\alpha_{i}}{\alpha_{i-1}} \qquad \qquad (1\leq i\leq d-1). \label{eq:ci_1}
\end{align}
\end{lemma}
\noindent {\it Proof:} First, eliminate $a^{*}_{0}$ and $c_{i}\alpha_{i-1}$ from (\ref{eq:genstailTTR2}) using (\ref{eq:a*0}) and (\ref{eq:cos_TTR2}) to obtain (\ref{eq:bi_1}). Substituting this back into (\ref{eq:cos_TTR2}), we obtain (\ref{eq:ci_1}). \hfill $\Box$ \\

\begin{lemma}
With reference to Assumption \ref{assum:A_bip_d>2}, further assume conditions {\rm (i)} and {\rm (iii)} from Note \ref{note:assumptions}. Suppose that in the diagram $\Delta$, vertex $0$ is adjacent to vertex $1$ and no other vertices. Then
\begin{align}
\theta_{0}\frac{\theta^{*}_{i-1}-\theta^{*}_{1}}{\theta^{*}_{i}-\theta^{*}_{0}} &\neq \theta_{1} \qquad \qquad (1\leq i\leq d-1), \label{eq:bine0} \\
\theta_{0}\frac{\theta^{*}_{i+1}-\theta^{*}_{1}}{\theta^{*}_{i}-\theta^{*}_{0}} &\neq \theta_{1} \qquad \qquad (1\leq i\leq d-1). \label{eq:cine0}
\end{align}
\end{lemma}
\noindent {\it Proof:} Without loss of generality, fix a feasible basis as in Assumption \ref{assum:fix_feas}. Suppose $\theta_{0}\frac{\theta^{*}_{i-1}-\theta^{*}_{1}}{\theta^{*}_{i}-\theta^{*}_{0}}=\theta_{1}$ for some integer $i$ such that $1\leq i\leq d-1$. Then $b_{i}=0$ by (\ref{eq:bi_1}). This is a contradiction, so (\ref{eq:bine0}) follows. Similarly, (\ref{eq:cine0}) is a consequence of (\ref{eq:ci_1}). \hfill $\Box$ \\

\section{Algebraic consequences of the tail condition} \label{sec:I}

\begin{definition} \label{def:dual_feasible}
\rm
With reference to Assumption \ref{assum:A_bip_d>2}, let $\{w_{i}\}_{i=0}^{d}$ denote a basis of $V$. We say that this basis is {\it dual feasible} whenever $w_{i}\in E_{i}V$ for $0\leq i\leq d$.
\end{definition}

\noindent With reference to Assumption \ref{assum:A_bip_d>2}, further assume conditions (ii) and (iii) from Note \ref{note:assumptions}. Let $\{w_{i}\}_{i=0}^{d}$ denote a dual feasible basis of $V$. By Proposition \ref{prop:EA*E} and Definition \ref{def:dual_feasible}, the matrices representing $A$ and $A^{*}$ with respect to $\{w_{i}\}_{i=0}^{d}$ are
\begin{equation} \label{eq:A_A*_mat_gen2}
A:\mbox{diag}\left(\theta_{0}, \theta_{1}, \ldots, \theta_{d}\right),
\qquad
A^{*}:
\left(
\begin{array}
{ c c}
A^{*}_{0,0} & A^{*}_{0,1} \\
A^{*}_{1,0} & A^{*}_{1,1} 
\end{array}
\right),
\end{equation}
where $A^{*}_{0,0}=\left(
\begin{array}
{ c c}
  a^{*}_{0} & b^{*}_{0} \\
  c^{*}_{1} & a^{*}_{1} 
\end{array}
\right)$, $A^{*}_{0,1}$ is the $2$ by $d-1$ matrix block with lower left corner entry $b^{*}_{1}$ and all other entries $0$, $A^{*}_{1,0}$ is the $d-1$ by $2$ matrix block with upper right corner entry $c^{*}_{2}$ and all other entries $0$, $A^{*}_{1,1}$ is a $d-1$ by $d-1$ matrix block with diagonal entries $a^{*}_{2}, \ldots, a^{*}_{d}$ and all other entries arbitrary, and the scalars $\{a^{*}_{i}\}_{i=0}^{d}$ are from Definition \ref{def:a*}. Because $\Delta$ is connected, each of the scalars $b^{*}_{0}$, $b^{*}_{1}$, $c^{*}_{1}$, and $c^{*}_{2}$ is nonzero by Definition \ref{def:delta}. For notational convenience, let $c^{*}_{0}=0$. Recall that $\{\theta_{i}\}_{i=0}^{d}$ are mutually distinct.

\medskip

\noindent Observe that the products $b^{*}_{0}c^{*}_{1}$ and $b^{*}_{1}c^{*}_{2}$ are independent of our choice of dual feasible basis.

\begin{proposition} \label{prop:tail_eq}
With reference to Assumptions \ref{assum:A_bip_d>2} and \ref{assum:fix_feas}, further assume conditions {\rm (ii)} and {\rm (iii)} from Note \ref{note:assumptions}. Then for $0\leq i\leq d$,
\begin{equation} \label{eq:gentailTTR2} 
\begin{split}
c_{i}(\theta^{*}_{i-1}-\theta^{*}_{0})(\theta^{*}_{i-1}-\theta^{*}_{1})\alpha_{i-1}+b_{i}(\theta^{*}_{i+1}-\theta^{*}_{0})(\theta^{*}_{i+1}-\theta^{*}_{1})\alpha_{i+1} \\
= \theta_{2}(\theta^{*}_{i}-\theta^{*}_{0})(\theta^{*}_{i}-\theta^{*}_{1})\alpha_{i}+(\theta_{2}-\theta_{1})\psi \theta^{*}_{i}\alpha_{i}+(\theta_{1}-\theta_{0})\psi a^{*}_{0}\alpha_{i}+(\theta_{2}-\theta_{0})\zeta\alpha_{i}, 
\end{split}
\end{equation}
where $\psi =\theta^{*}_{1}+\theta^{*}_{0}-a^{*}_{1}-a^{*}_{0}$ and $\zeta =a^{*}_{0}a^{*}_{1}-b^{*}_{0}c^{*}_{1}-\theta^{*}_{0}\theta^{*}_{1}$.
\end{proposition}
\noindent {\it Proof:} Fix a dual feasible basis $\{v^{*}_{i}\}_{i=0}^{d}$ of $V$ such that $v^{*}_{0}=\sum_{i=0}^{d}\alpha_{i}v_{i}$. Let $w=(A^{*}-a^{*}_{0}I)v^{*}_{0}$ and note that $w=\sum_{i=0}^{d}\alpha'_{i}v_{i}$ where $\alpha'_{i}=(\theta^{*}_{i}-a^{*}_{0})\alpha_{i}$ for $0\leq i\leq d$. By Lemma \ref{lem:A*-a*I}, $w\in E_{1}V$. Because  $(E_{0}, E_{1})$ is a tail in $\Delta$ connected to vertex $2$, it follows that $E_{i}A^{*}E_{1}=0$ for $3\leq i\leq d$. So, $A^{*}E_{1}V\subseteq E_{0}V+E_{1}V+E_{2}V$. With respect to the $\{v^{*}_{i}\}_{i=0}^{d}$ basis, $A$ and $A^{*}$ are represented by the matrices in (\ref{eq:A_A*_mat_gen2}). Using matrix multiplication, we calculate $w$ with respect to this basis,
\begin{align}
w &= \left(A^{*}-a^{*}_{0}I\right)v^{*}_{0} \notag \\
  &= a^{*}_{0}v^{*}_{0}+c^{*}_{1}v^{*}_{1} - a^{*}_{0}v^{*}_{0} \notag \\
  &= c^{*}_{1}v^{*}_{1}. \notag
\end{align}
Accordingly, 
\begin{equation}
A^{*}w = b^{*}_{0}c^{*}_{1}v^{*}_{0}+a^{*}_{1}c^{*}_{1}v^{*}_{1}+c^{*}_{2}c^{*}_{1}v^{*}_{2}. \notag
\end{equation}
Therefore, $A^{*}w-a^{*}_{1}w-b^{*}_{0}c^{*}_{1}v^{*}_{0}\in E_{2}V$. Substituting $w=(A^{*}-a^{*}_{0}I)v^{*}_{0}$, we obtain $(A^{*}-a^{*}_{1}I)(A^{*}-a^{*}_{0}I)v^{*}_{0}-b^{*}_{0}c^{*}_{1}v^{*}_{0}\in E_{2}V$. Therefore,
\begin{equation}
A((A^{*}-a^{*}_{1}I)(A^{*}-a^{*}_{0}I)-b^{*}_{0}c^{*}_{1})v^{*}_{0}=\theta_{2}((A^{*}-a^{*}_{1}I)(A^{*}-a^{*}_{0}I)-b^{*}_{0}c^{*}_{1})v^{*}_{0}. \notag
\end{equation}
Recall that with respect to the basis $\{v_{i}\}_{i=0}^{d}$, $A$ and $A^{*}$ are represented by the matrices from (\ref{eq:A_A*_mat_gen}), and $v^{*}_{0}$ is represented by $(\alpha_{0},\alpha_{1},\ldots ,\alpha_{d})^{t}$. The result follows by matrix multiplication. \hfill $\Box$ \\

\noindent We now proceed to derive expressions for the scalars $\psi$ and $\zeta$ from (\ref{eq:gentailTTR2}).

\begin{lemma} 
With reference to Assumption \ref{assum:A_bip_d>2}, further assume conditions {\rm (i)--(iii)} from Note \ref{note:assumptions}. Then the scalars $\psi$ and $\zeta$ from Proposition \ref{prop:tail_eq} satisfy 
\begin{align}
\psi &= \frac{\theta_{1}(\theta^{*}_{d-1}-\theta^{*}_{0})-\theta_{2}(\theta^{*}_{d}-\theta^{*}_{1})}{\theta_{2}-\theta_{1}}, \label{eq:a*0+a*1} \\
\zeta &= \frac{\theta_{2}\theta^{*}_{0}-\theta_{0}\theta^{*}_{1}}{\theta_{2}-\theta_{0}}\frac{\theta_{2}(\theta^{*}_{d}-\theta^{*}_{1})-\theta_{1}(\theta^{*}_{d-1}-\theta^{*}_{0})}{\theta_{2}-\theta_{1}}. \label{eq:a*0a*1-b*0c*1}
\end{align}
\end{lemma}
\noindent {\it Proof:} Without loss of generality, fix a feasible basis as in Assumption \ref{assum:fix_feas}. Consider the two equations obtained from taking (\ref{eq:gentailTTR2}) at $i=0$ and $i=d$. Eliminate $b_{0}$ using (\ref{eq:b0alpha1}), $c_{d}$ using (\ref{eq:cdalphad-1}), and $a^{*}_{0}$ using (\ref{eq:a*0}). Invoking the fact that both $\alpha_{0}$ and $\alpha_{d}$ are nonzero by Note \ref{note:assumptions}(i), we obtain the following system of equations
\begin{align}
(\theta_{2}\theta^{*}_{0}-\theta_{0}\theta^{*}_{1})\psi +(\theta_{2}-\theta_{0})\zeta &= 0, \label{eq:psizeta1} \\
(\theta_{2}\theta^{*}_{d}-\theta_{0}\theta^{*}_{d-1})\psi +(\theta_{2}-\theta_{0})\zeta &= \theta_{0}(\theta^{*}_{d-1}-\theta^{*}_{0})(\theta^{*}_{d-1}-\theta^{*}_{1})-\theta_{2}(\theta^{*}_{d}-\theta^{*}_{0})(\theta^{*}_{d}-\theta^{*}_{1}). \label{eq:psizeta2}
\end{align}
Take the difference between (\ref{eq:psizeta1}) and (\ref{eq:psizeta2}). Solve the resulting equation for $\psi$. In the expression for $\psi$, eliminate $\theta_{0}$ using (\ref{eq:theta01}) and simplify to obtain (\ref{eq:a*0+a*1}).
Using this and (\ref{eq:psizeta2}), we obtain (\ref{eq:a*0a*1-b*0c*1}). \hfill $\Box$ \\

\begin{proposition} \label{prop:theta012}
With reference to Assumption \ref{assum:A_bip_d>2}, further assume conditions {\rm (i)--(iii)} from Note \ref{note:assumptions}. Then for $1\leq i\leq d-1$,
\begin{equation}
\begin{split}
\theta_{0}(\theta^{*}_{i+1}-\theta^{*}_{1})(\theta^{*}_{i-1}-\theta^{*}_{1})+\theta_{1}(\theta^{*}_{i}-\theta^{*}_{0})(\theta^{*}_{d-1}-\theta^{*}_{i+1}-\theta^{*}_{i-1}+\theta^{*}_{1}) \\
-\theta_{2}(\theta^{*}_{i}-\theta^{*}_{0})(\theta^{*}_{d}-\theta^{*}_{i})=0. \label{eq:theta012_gen}
\end{split}
\end{equation}
\end{proposition}
\noindent {\it Proof:} Without loss of generality, fix a feasible basis as in Assumption \ref{assum:fix_feas}. Eliminate $a^{*}_{0}$, $\psi$, $\zeta$, $b_{i}$, and $c_{i}$ from (\ref{eq:gentailTTR2}) using (\ref{eq:a*0}), (\ref{eq:a*0+a*1}), (\ref{eq:a*0a*1-b*0c*1}), (\ref{eq:bi_1}), and (\ref{eq:ci_1}), respectively. The result follows after simplifying and invoking the fact that $\alpha_{i}\neq 0$ ($0\leq i\leq d$) by Note \ref{note:assumptions}(i). \hfill $\Box$ \\

\section{The three-term recurrence}

We begin by recalling the following fact from arithmetic, whose proof is left as an exercise for the reader.

\begin{lemma} \label{lem:frac_arith}
Let $m$, $n$, $r$, and $s$ denote elements of $\fld$ such that $n$, $s$, and $n+s$ are nonzero and $m/n=r/s$. Then this common value equals $(m+r)/(n+s)$.
\end{lemma}

\noindent We now utilize Lemma \ref{lem:frac_arith} to establish some algebraic facts for later use. With reference to Assumption \ref{assum:A_bip_d>2}, our goal is to show that under assumptions (i)--(iii) from Note \ref{note:assumptions},
\begin{equation}
\frac{\theta^{*}_{j}-\theta^{*}_{j-3}}{\theta^{*}_{j-1}-\theta^{*}_{j-2}} \notag
\end{equation}
is independent of $j$ for $3\leq j\leq d$.

\begin{proposition} \label{prop:alg_facts1}
With reference to Assumption \ref{assum:A_bip_d>2}, further assume condition {\rm (iii)} from Note \ref{note:assumptions}. Fix an integer $i$ such that $3\leq i\leq d$. Then the following {\rm (i)}, {\rm (ii)} are equivalent.
\begin{enumerate}
\item[\rm (i)] The expression
\begin{equation} \label{eq:frac_recur1}
\frac{\theta^{*}_{j}-\theta^{*}_{j-3}}{\theta^{*}_{j-1}-\theta^{*}_{j-2}}
\end{equation}
is independent of $j$ for $3\leq j\leq i$.
\item[\rm (ii)] The expression
\begin{equation} \label{eq:frac_recur2}
\frac{\theta^{*}_{j}+\theta^{*}_{j-1}+\theta^{*}_{j-2}-\theta^{*}_{h}-\theta^{*}_{h-1}-\theta^{*}_{h-2}}{\theta^{*}_{j-1}-\theta^{*}_{h-1}}
\end{equation}
is independent of $j$ and $h$ for $2\leq h< j\leq i$.
\end{enumerate}
Suppose conditions {\rm (i)} and {\rm (ii)} hold. Then the common values of {\rm (\ref{eq:frac_recur1})} and {\rm (\ref{eq:frac_recur2})} are equal.
\end{proposition}
\noindent {\it Proof:} (i) $\Rightarrow$ (ii). Note that for any $j$ and $h$ such that $3\leq j\leq i$ and $2\leq h\leq j-1$,
\begin{equation} \label{eq:frac_chain1}
\frac{\theta^{*}_{j}-\theta^{*}_{j-3}}{\theta^{*}_{j-1}-\theta^{*}_{j-2}}=\frac{\theta^{*}_{j-1}-\theta^{*}_{j-4}}{\theta^{*}_{j-2}-\theta^{*}_{j-3}}=\cdots =\frac{\theta^{*}_{h+2}-\theta^{*}_{h-1}}{\theta^{*}_{h+1}-\theta^{*}_{h}}=\frac{\theta^{*}_{h+1}-\theta^{*}_{h-2}}{\theta^{*}_{h}-\theta^{*}_{h-1}}.
\end{equation}
A consequence of Lemma \ref{lem:frac_arith} is that all of the equal fractions in (\ref{eq:frac_chain1}) also equal the fraction whose numerator is the sum of all the numerators appearing in (\ref{eq:frac_chain1}) and whose denominator is the sum of all the denominators appearing in (\ref{eq:frac_chain1}). Note that both of these sums are telescoping and (ii) follows.

\medskip

\noindent (ii) $\Rightarrow$ (i). For any $j$ such that $3\leq j\leq i$, let $h=j-1$. Then (i) follows immediately.

\medskip

\noindent Now suppose (i) and (ii) hold. It follows from the above arguments that the common values of (\ref{eq:frac_recur1}) and (\ref{eq:frac_recur2}) are equal. \hfill $\Box$ \\

\begin{lemma} \label{lem:alg_facts2}
With reference to Assumption \ref{assum:A_bip_d>2}, further assume condition {\rm (iii)} from Note \ref{note:assumptions}. Fix an integer $i$ such that $3\leq i\leq d$ and assume that the equivalent conditions {\rm (i)} and {\rm (ii)} from Proposition \ref{prop:alg_facts1} hold. Then  
\begin{equation} \label{eq:extra_id}
\frac{\theta^{*}_{i-2}-\theta^{*}_{2}}{\theta^{*}_{i-1}-\theta^{*}_{1}}=\frac{\theta^{*}_{i-1}-\theta^{*}_{3}}{\theta^{*}_{i}-\theta^{*}_{2}}.
\end{equation}
\end{lemma}
\noindent {\it Proof:} First note that for $i=3$ and $i=4$, (\ref{eq:extra_id}) is trivial. If $i=5$, then (\ref{eq:extra_id}) is true by assumption. So, we assume without loss of generality that $i\geq 6$. We consider two cases, depending on whether $i$ is even or $i$ is odd.

\medskip

\noindent {\it Case 1: $i=2n$ is even.} By assumption, note that 
\begin{equation}
\frac{\theta^{*}_{n+2}+\theta^{*}_{n+1}+\theta^{*}_{n}-\theta^{*}_{n}-\theta^{*}_{n-1}-\theta^{*}_{n-2}}{\theta^{*}_{n+1}-\theta^{*}_{n-1}} = \frac{\theta^{*}_{n+3}+\theta^{*}_{n+2}+\theta^{*}_{n+1}-\theta^{*}_{n+1}-\theta^{*}_{n}-\theta^{*}_{n-1}}{\theta^{*}_{n+2}-\theta^{*}_{n}}, \notag
\end{equation}
which further simplifies to
\begin{equation}
\frac{\theta^{*}_{n+2}-\theta^{*}_{n-2}}{\theta^{*}_{n+1}-\theta^{*}_{n-1}}=\frac{\theta^{*}_{n+3}-\theta^{*}_{n-1}}{\theta^{*}_{n+2}-\theta^{*}_{n}}. \notag
\end{equation}
Taking the reciprocal of both sides, we obtain
\begin{equation} \label{eq:frac1e}
\frac{\theta^{*}_{n+1}-\theta^{*}_{n-1}}{\theta^{*}_{n+2}-\theta^{*}_{n-2}}=\frac{\theta^{*}_{n+2}-\theta^{*}_{n}}{\theta^{*}_{n+3}-\theta^{*}_{n-1}}.
\end{equation}
If $i=6$, then (\ref{eq:frac1e}) is (\ref{eq:extra_id}) and we are done. Otherwise, it also follows by assumption that 
\begin{equation} \label{eq:pre_frac2e}
\frac{\theta^{*}_{n+3}+\theta^{*}_{n+2}+\theta^{*}_{n+1}-\theta^{*}_{n-1}-\theta^{*}_{n-2}-\theta^{*}_{n-3}}{\theta^{*}_{n+2}-\theta^{*}_{n-2}} = \frac{\theta^{*}_{n+4}+\theta^{*}_{n+3}+\theta^{*}_{n+2}-\theta^{*}_{n}-\theta^{*}_{n-1}-\theta^{*}_{n-2}}{\theta^{*}_{n+3}-\theta^{*}_{n-1}}.
\end{equation}
Using (\ref{eq:frac1e}) to simplify (\ref{eq:pre_frac2e}), we obtain
\begin{equation}
\frac{\theta^{*}_{n+3}-\theta^{*}_{n-3}}{\theta^{*}_{n+2}-\theta^{*}_{n-2}}=\frac{\theta^{*}_{n+4}-\theta^{*}_{n-2}}{\theta^{*}_{n+3}-\theta^{*}_{n-1}}. \notag
\end{equation}
Again, taking the reciprocal of both sides, we obtain
\begin{equation}
\frac{\theta^{*}_{n+2}-\theta^{*}_{n-2}}{\theta^{*}_{n+3}-\theta^{*}_{n-3}}=\frac{\theta^{*}_{n+3}-\theta^{*}_{n-1}}{\theta^{*}_{n+4}-\theta^{*}_{n-2}}. \notag
\end{equation}
If $i=8$, then we are done. Otherwise, we continue in the above fashion until we eventually obtain (\ref{eq:extra_id}).

\medskip

\noindent {\it Case 2: $i=2n+1$ is odd.} By assumption, it follows that 
\begin{equation} \label{eq:frac0o}
\frac{\theta^{*}_{n+1}-\theta^{*}_{n}}{\theta^{*}_{n+2}-\theta^{*}_{n-1}}=\frac{\theta^{*}_{n+2}-\theta^{*}_{n+1}}{\theta^{*}_{n+3}-\theta^{*}_{n}}
\end{equation}
and
\begin{equation} \label{eq:pre_frac1o}
\frac{\theta^{*}_{n+3}+\theta^{*}_{n+2}+\theta^{*}_{n+1}-\theta^{*}_{n}-\theta^{*}_{n-1}-\theta^{*}_{n-2}}{\theta^{*}_{n+2}-\theta^{*}_{n-1}} = \frac{\theta^{*}_{n+4}+\theta^{*}_{n+3}+\theta^{*}_{n+2}-\theta^{*}_{n+1}-\theta^{*}_{n}-\theta^{*}_{n-1}}{\theta^{*}_{n+3}-\theta^{*}_{n}}.
\end{equation}
Using (\ref{eq:frac0o}) to simplify (\ref{eq:pre_frac1o}), we obtain
\begin{equation}
\frac{\theta^{*}_{n+3}-\theta^{*}_{n-2}}{\theta^{*}_{n+2}-\theta^{*}_{n-1}}=\frac{\theta^{*}_{n+4}-\theta^{*}_{n-1}}{\theta^{*}_{n+3}-\theta^{*}_{n}}. \notag
\end{equation}
Taking the reciprocal of both sides, we obtain
\begin{equation} \label{eq:frac1o}
\frac{\theta^{*}_{n+2}-\theta^{*}_{n-1}}{\theta^{*}_{n+3}-\theta^{*}_{n-2}}=\frac{\theta^{*}_{n+3}-\theta^{*}_{n}}{\theta^{*}_{n+4}-\theta^{*}_{n-1}}.
\end{equation}
If $i=7$, then (\ref{eq:frac1o}) is (\ref{eq:extra_id}) and we are done. Otherwise, it also follows by assumption that 
\begin{equation} \label{eq:pre_frac2o}
\frac{\theta^{*}_{n+4}+\theta^{*}_{n+3}+\theta^{*}_{n+2}-\theta^{*}_{n-1}-\theta^{*}_{n-2}-\theta^{*}_{n-3}}{\theta^{*}_{n+3}-\theta^{*}_{n-2}} = \frac{\theta^{*}_{n+5}+\theta^{*}_{n+4}+\theta^{*}_{n+3}-\theta^{*}_{n}-\theta^{*}_{n-1}-\theta^{*}_{n-2}}{\theta^{*}_{n+4}-\theta^{*}_{n-1}}.
\end{equation}
Using (\ref{eq:frac1o}) to simplify (\ref{eq:pre_frac2o}), we obtain
\begin{equation}
\frac{\theta^{*}_{n+4}-\theta^{*}_{n-3}}{\theta^{*}_{n+3}-\theta^{*}_{n-2}}=\frac{\theta^{*}_{n+5}-\theta^{*}_{n-2}}{\theta^{*}_{n+4}-\theta^{*}_{n-1}}. \notag
\end{equation}
Again, taking the reciprocal of both sides, we obtain
\begin{equation}
\frac{\theta^{*}_{n+3}-\theta^{*}_{n-2}}{\theta^{*}_{n+4}-\theta^{*}_{n-3}}=\frac{\theta^{*}_{n+4}-\theta^{*}_{n-1}}{\theta^{*}_{n+5}-\theta^{*}_{n-2}}. \notag
\end{equation}
If $i=9$, then we are done. Otherwise, we continue in the above fashion until we eventually obtain (\ref{eq:extra_id}). \hfill $\Box$ \\

\begin{proposition} \label{prop:recur}
With reference to Assumption \ref{assum:A_bip_d>2}, further assume conditions {\rm (i)--(iii)} from Note \ref{note:assumptions}. Then
\begin{equation} \label{eq:TTR_frac}
\frac{\theta^{*}_{j}-\theta^{*}_{j-3}}{\theta^{*}_{j-1}-\theta^{*}_{j-2}}
\end{equation}
is independent of $j$ for $3\leq j\leq d$.
\end{proposition}
\noindent {\it Proof:} First note that for $d=3$, there is nothing to prove, so we may assume without loss of generality that $d\geq 4$. We begin by evaluating (\ref{eq:theta012_gen}) at $i=2$ and solving for $\theta_{1}\theta^{*}_{d-1}$ to obtain
\begin{equation} \label{eq:elimd}
\theta_{1}\theta^{*}_{d-1}=\theta_{2}(\theta^{*}_{d}-\theta^{*}_{2})+\theta_{1}\theta^{*}_{3}.
\end{equation}
We now use (\ref{eq:elimd}) to eliminate $\theta_{1}\theta^{*}_{d-1}$ from (\ref{eq:theta012_gen}) for $i=1$ and $3\leq i\leq d-1$. For $i=1$, after simplifying we obtain
\begin{equation}
\theta_{2}=\theta_{1}\left(\frac{\theta^{*}_{3}-\theta^{*}_{0}}{\theta^{*}_{2}-\theta^{*}_{1}}-1\right)-\theta_{0}. \label{eq:theta012_1}
\end{equation}
For $3\leq i\leq d-1$, we obtain
\begin{equation}
\theta_{0}\frac{\theta^{*}_{i-1}-\theta^{*}_{1}}{\theta^{*}_{i}-\theta^{*}_{0}}-\theta_{1}\frac{\theta^{*}_{i+1}+\theta^{*}_{i-1}-\theta^{*}_{3}-\theta^{*}_{1}}{\theta^{*}_{i+1}-\theta^{*}_{1}}+\theta_{2}\frac{\theta^{*}_{i}-\theta^{*}_{2}}{\theta^{*}_{i+1}-\theta^{*}_{1}}=0. \label{eq:theta012_genalt}
\end{equation}
Eliminating $\theta_{2}$ from (\ref{eq:theta012_genalt}) using (\ref{eq:theta012_1}), we obtain
\begin{equation} \label{eq:theta01_genalti}
\theta_{0}\frac{\theta^{*}_{i-1}-\theta^{*}_{1}}{\theta^{*}_{i}-\theta^{*}_{0}}\left(\frac{\theta^{*}_{i+1}-\theta^{*}_{1}}{\theta^{*}_{i}-\theta^{*}_{2}}-\frac{\theta^{*}_{i}-\theta^{*}_{0}}{\theta^{*}_{i-1}-\theta^{*}_{1}}\right)=\theta_{1}T^{*}_{i}, \qquad \qquad (3\leq i\leq d-1)
\end{equation}
where 
\begin{equation}
T^{*}_{i}=\frac{\theta^{*}_{i+1}+\theta^{*}_{i}+\theta^{*}_{i-1}-\theta^{*}_{3}-\theta^{*}_{2}-\theta^{*}_{1}}{\theta^{*}_{i}-\theta^{*}_{2}}-\frac{\theta^{*}_{3}-\theta^{*}_{0}}{\theta^{*}_{2}-\theta^{*}_{1}}. \notag
\end{equation}
Assume by way of a contradiction that (\ref{eq:TTR_frac}) is not independent of $j$ for $3\leq j\leq d$. Let $k$ denote the smallest integer ($4\leq k\leq d$) such that (\ref{eq:TTR_frac}) at $j=k$ does not equal (\ref{eq:TTR_frac}) at $j=k-1$. By construction, (\ref{eq:TTR_frac}) is independent of $j$ for $3\leq j< k$. Therefore, Proposition~\ref{prop:alg_facts1}(i) holds with $i=k-1$. Consequently, 
\begin{equation}
\frac{\theta^{*}_{3}-\theta^{*}_{0}}{\theta^{*}_{2}-\theta^{*}_{1}}=\frac{\theta^{*}_{k-1}+\theta^{*}_{k-2}+\theta^{*}_{k-3}-\theta^{*}_{2}-\theta^{*}_{1}-\theta^{*}_{0}}{\theta^{*}_{k-2}-\theta^{*}_{1}}. \notag
\end{equation}
Using this and Lemma \ref{lem:alg_facts2} with $i=k-1$, 
\begin{align}
T^{*}_{k-1} &= \frac{\theta^{*}_{k}-\theta^{*}_{1}}{\theta^{*}_{k-1}-\theta^{*}_{2}}-\frac{\theta^{*}_{k-1}-\theta^{*}_{0}}{\theta^{*}_{k-2}-\theta^{*}_{1}}+\frac{\theta^{*}_{k-2}-\theta^{*}_{3}}{\theta^{*}_{k-1}-\theta^{*}_{2}}-\frac{\theta^{*}_{k-3}-\theta^{*}_{2}}{\theta^{*}_{k-2}-\theta^{*}_{1}}, \notag \\
            &= \frac{\theta^{*}_{k}-\theta^{*}_{1}}{\theta^{*}_{k-1}-\theta^{*}_{2}}-\frac{\theta^{*}_{k-1}-\theta^{*}_{0}}{\theta^{*}_{k-2}-\theta^{*}_{1}}. \notag
\end{align}
So by (\ref{eq:theta01_genalti}) with $i=k-1$, $T_{k-1}^{*}$ times
\begin{equation} \label{eq:bi_num}
\theta_{0}\frac{\theta^{*}_{k-2}-\theta^{*}_{1}}{\theta^{*}_{k-1}-\theta^{*}_{0}}-\theta_{1}
\end{equation}
equals $0$. Therefore, either $T_{k-1}^{*}=0$ or the expression (\ref{eq:bi_num}) equals $0$. By (\ref{eq:bine0}), the former must be true. Combining this with the fact that (\ref{eq:TTR_frac}) is independent of $j$ for $3\leq j< k$, 
\begin{equation}
\frac{\theta^{*}_{3}-\theta^{*}_{0}}{\theta^{*}_{2}-\theta^{*}_{1}}=\frac{\theta^{*}_{4}-\theta^{*}_{1}}{\theta^{*}_{3}-\theta^{*}_{2}}=\cdots=\frac{\theta^{*}_{k-1}-\theta^{*}_{k-4}}{\theta^{*}_{k-2}-\theta^{*}_{k-3}}=\frac{\theta^{*}_{k}+\theta^{*}_{k-1}+\theta^{*}_{k-2}-\theta^{*}_{3}-\theta^{*}_{2}-\theta^{*}_{1}}{\theta^{*}_{k-1}-\theta^{*}_{2}}. \notag
\end{equation}
By Lemma \ref{lem:frac_arith},
\begin{equation}
\frac{\theta^{*}_{3}-\theta^{*}_{0}}{\theta^{*}_{2}-\theta^{*}_{1}}=\frac{\theta^{*}_{4}-\theta^{*}_{1}}{\theta^{*}_{3}-\theta^{*}_{2}}=\cdots=\frac{\theta^{*}_{k-1}-\theta^{*}_{k-4}}{\theta^{*}_{k-2}-\theta^{*}_{k-3}}=\frac{\theta^{*}_{k}-\theta^{*}_{k-3}}{\theta^{*}_{k-1}-\theta^{*}_{k-2}}. \notag
\end{equation}
In other words, (\ref{eq:TTR_frac}) is independent of $j$ for $3\leq j\leq k$. This is a contradiction, so our result follows. \hfill $\Box$ \\

\begin{corollary} \label{cor:TTR}
With reference to Assumption \ref{assum:A_bip_d>2}, further assume conditions {\rm (i)--(iii)} from Note \ref{note:assumptions}. Then there exists $\beta \in \fld$ such that $\theta^{*}_{i-1}-\beta \theta^{*}_{i}+\theta^{*}_{i+1}$ is independent of $i$ for $1\leq i \leq d-1$.
\end{corollary}
\noindent {\it Proof:} This is an immediate consequence of Proposition \ref{prop:recur} combined with \cite[Lemma 8.3]{T:Leonard}. \hfill $\Box$ \\

\section{The main theorem}

The following is our main theorem.

\begin{theorem} \label{thm:main}
With reference to Assumption \ref{assum:A_bip}, let $(E,F)$ denote an ordered pair of distinct primitive idempotents for $A$. Then this pair is $Q$-polynomial if and only if the following {\rm (i)--(iii)} hold.
\begin{enumerate}
\item[\rm (i)] The primitive idempotent $E$ is normalizing.
\item[\rm (ii)] $(E,F)$ is a tail.
\item[\rm (iii)] $\{\theta^{*}_{i}\}_{i=0}^{d}$ are mutually distinct.
\end{enumerate}
\end{theorem}
\noindent {\it Proof:} First, assume that $(E,F)$ is $Q$-polynomial. Fix the $Q$-polynomial ordering $\{E_{i}\}_{i=0}^{d}$ such that $E=E_{0}$. Condition (i) follows from \cite[Lemma 10.7]{T:madrid}. Condition (ii) follows from Lemma \ref{lem:Qpolytail}. Condition (iii) follows from the last paragraph of Section \ref{sec:LS}.

\medskip

\noindent Conversely, assume conditions (i)--(iii). If $d<3$, then $\Delta$ is clearly a path. Otherwise, there exists $\beta \in \fld$ such that $\theta^{*}_{i-1}-\beta \theta^{*}_{i}+\theta^{*}_{i+1}$ is independent of $i$ for $1\leq i \leq d-1$ by Corollary~\ref{cor:TTR}. By Theorem \ref{thm:tail}, it follows that $(E,F)$ is $Q$-polynomial. \hfill $\Box$ \\

\section{Bipartite Leonard systems}

We conclude with an observation regarding bipartite Leonard systems. We first establish a necessary lemma.

\begin{lemma} \label{lem:theta12}
With reference to Assumption \ref{assum:A_bip_d>2}, further assume conditions {\rm (i)--(iii)} from Note \ref{note:assumptions}. Then
\begin{equation} \label{eq:theta12}
\theta_{1}(\theta^{*}_{d-1}-\theta^{*}_{3})=\theta_{2}(\theta^{*}_{d}-\theta^{*}_{2}).
\end{equation}
\end{lemma}
\noindent {\it Proof:} This is a direct consequence of (\ref{eq:theta012_gen}) at $i=2$. \hfill $\Box$ \\

\begin{lemma} \label{lem:theta_sum}
With reference to Assumption \ref{assum:A_bip}, if $(A; \{E_{i}\}^{d}_{i=0}; A^{*}; \{E^{*}_{i}\}^{d}_{i=0})$ is a Leonard system then
\begin{equation} \label{eq:theta_opp_sym}
\theta_{i}+\theta_{d-i}=0 \qquad \qquad (0\leq i\leq d).
\end{equation}
\end{lemma}
\noindent {\it Proof:} Combine Assumption \ref{assum:A_bip} and Definition \ref{def:bipartite} with \cite[Lemma 10.2]{NT:affine}. \hfill $\Box$ \\

\begin{corollary}
With reference to Assumption \ref{assum:A_bip}, if $(A; \{E_{i}\}^{d}_{i=0}; A^{*}; \{E^{*}_{i}\}^{d}_{i=0})$ is a Leonard system then $\theta_{i}=0$ if and only if $d$ is even and $i=d/2$.
\end{corollary}
\noindent {\it Proof:} This follows immediately from Lemma \ref{lem:theta_sum} and the fact that $\{\theta_{i}\}_{i=0}^{d}$ are mutually distinct. \hfill $\Box$ \\

\begin{proposition} \label{prop:theta_i,i+1}
With reference to Assumption \ref{assum:A_bip}, if $(A; \{E_{i}\}^{d}_{i=0}; A^{*}; \{E^{*}_{i}\}^{d}_{i=0})$ is a Leonard system then
\begin{equation} \label{eq:theta_i,i+1}
\theta_{i}(\theta^{*}_{d-i-1}-\theta^{*}_{i+1}) = \theta_{i+1}(\theta^{*}_{d-i}-\theta^{*}_{i}) \qquad \qquad (0\leq i\leq d-1).
\end{equation}
\end{proposition}
\noindent {\it Proof:} Suppose $d=1$ or $d=2$. Then (\ref{eq:theta_i,i+1}) follows from (\ref{eq:theta_opp_sym}) in both cases. So, assume $d\geq 3$. We proceed by induction on $i$. We first establish two base cases that we will use in a two-step induction hypothesis. By Lemma \ref{lem:lspath}, $(E_{0}, E_{1})$ is a tail in $\Delta$ and vertex $1$ is connected to vertex $2$ in $\Delta$. By Theorem \ref{thm:main}, conditions (i)--(iii) from Note \ref{note:assumptions} hold. Therefore, Proposition \ref{prop:theta012} holds. For $i=0$, (\ref{eq:theta_i,i+1}) is just (\ref{eq:theta01}). For $i=1$, (\ref{eq:theta_i,i+1}) follows from (\ref{eq:theta12}) combined with \cite[Lemma 9.4]{T:Leonard}. For $2\leq i\leq d-1$, recall that the expressions in (\ref{eq:thetarecur}) are equal, so 
\begin{equation}
\frac{\theta_{i+1}-\theta_{i-2}}{\theta_{i}-\theta_{i-1}}=\frac{\theta^{*}_{i+1}-\theta^{*}_{i-2}}{\theta^{*}_{i}-\theta^{*}_{i-1}}. \notag
\end{equation}
From this,
\begin{equation}
\theta_{i+1}-\theta_{i-2}=\frac{\theta^{*}_{i+1}-\theta^{*}_{i-2}}{\theta^{*}_{i}-\theta^{*}_{i-1}}(\theta_{i}-\theta_{i-1}). \notag
\end{equation}
We now multiply each side by $\theta^{*}_{d-i+1}-\theta^{*}_{i-1}$ and invoke the induction hypothesis of $\theta_{i-2}(\theta^{*}_{d-i+1}-\theta^{*}_{i-1})=\theta_{i-1}(\theta^{*}_{d-i+2}-\theta^{*}_{i-2})$ to obtain
\begin{equation}
\theta_{i+1}(\theta^{*}_{d-i+1}-\theta^{*}_{i-1})=\frac{\theta^{*}_{i+1}-\theta^{*}_{i-2}}{\theta^{*}_{i}-\theta^{*}_{i-1}}(\theta_{i}-\theta_{i-1})(\theta^{*}_{d-i+1}-\theta^{*}_{i-1})+\theta_{i-1}(\theta^{*}_{d-i+2}-\theta^{*}_{i-2}). \notag
\end{equation}
Now multiply each side by $\theta^{*}_{d-i}-\theta^{*}_{i}$ and invoke the induction hypothesis of $\theta_{i-1}(\theta^{*}_{d-i}-\theta^{*}_{i})=\theta_{i}(\theta^{*}_{d-i+1}-\theta^{*}_{i-1})$. After simplifying, we obtain
\begin{equation}
\theta_{i+1}(\theta^{*}_{d-i}-\theta^{*}_{i})=\theta_{i}\left(-\frac{\theta^{*}_{i+1}-\theta^{*}_{i-2}}{\theta^{*}_{i}-\theta^{*}_{i-1}}(\theta^{*}_{d-i+1}-\theta^{*}_{d-i})-\theta^{*}_{i+1}+\theta^{*}_{d-i+2}\right). \notag
\end{equation}
The result follows by using 
\begin{equation}
\frac{\theta^{*}_{i+1}-\theta^{*}_{i-2}}{\theta^{*}_{i}-\theta^{*}_{i-1}}=\frac{\theta^{*}_{d-i+2}-\theta^{*}_{d-i-1}}{\theta^{*}_{d-i+1}-\theta^{*}_{d-i}} \notag
\end{equation}
and simplifying. \hfill $\Box$ \\

\begin{note}
\rm
Given a bipartite Leonard system $(A; \{E_{i}\}^{d}_{i=0}; A^{*}; \{E^{*}_{i}\}^{d}_{i=0})$, (\ref{eq:theta_i,i+1}) gives a way to recursively calculate $\{\theta_{i}\}_{i=0}^{d}$. Specifically, choose $j$ such that $0\leq j\leq d$. Provided $\theta_{j}\neq 0$, the sequence $\{\theta_{i}\}_{i=0}^{d}$ can be computed using $\{\theta^{*}_{i}\}_{i=0}^{d}$, $\theta_{j}$, Lemma \ref{lem:theta_sum}, and Proposition \ref{prop:theta_i,i+1}.
\end{note}

\section{Acknowledgment}

This paper was written while the author was a graduate student at the University of Wisconsin-Madison. The author thanks his advisor, Paul Terwilliger, for offering many valuable ideas and suggestions.

\bigskip

\noindent Edward Hanson \hfil\break
\noindent Department of Mathematics \hfil\break
\noindent University of Wisconsin \hfil\break
\noindent 480 Lincoln Drive \hfil\break
\noindent Madison, WI 53706-1388 USA \hfil\break
\noindent email: {\tt hanson@math.wisc.edu }\hfil\break

\end{document}